\documentclass[11pt]{article}

\usepackage{graphics}
\usepackage{amsfonts}
\usepackage{amssymb}
\usepackage{amsmath}
\usepackage{amsthm}
\usepackage{url}
\usepackage{color}
\usepackage{multirow}
\usepackage{enumerate}
\usepackage{hyperref}

\newtheorem{thm}{Theorem}

\newtheorem{prop}[thm]{Proposition}

\newtheorem{conj}[thm]{Conjecture}

\theoremstyle{definition}

\theoremstyle{definition}

\theoremstyle{definition}
\newtheorem{example}{Example}

\newcommand{\R}{\mathbb{R}}

\DeclareMathOperator{\tr}{tr}
\DeclareMathOperator{\diag}{diag}
\DeclareMathOperator{\Diag}{Diag}

\DeclareMathOperator{\rk}{rank}

\renewcommand{\binom}[2]{\genfrac(){0pt}{}{#1}{#2}}

\title{\Large{New bounds for the max-$k$-cut and chromatic number of a graph}}
\author{E.R. van Dam\thanks{Department of Econometrics and OR, Tilburg
University, The Netherlands. {\tt edwin.vandam@uvt.nl} }
 \and {R. Sotirov}\thanks{Department of Econometrics and OR, Tilburg
University, The Netherlands. {\tt r.sotirov@uvt.nl} }}

\date{}

\begin{document}
\maketitle

\begin{abstract}
We consider several semidefinite programming relaxations for the max-$k$-cut problem, with increasing complexity.
The optimal solution of the weakest presented semidefinite programming relaxation has a closed form
expression that includes the largest Laplacian eigenvalue of the graph under consideration.
This is the first known eigenvalue bound for the max-$k$-cut when $k>2$ that is applicable to any graph.
This bound is exploited to derive a new eigenvalue bound on the chromatic number of a graph.
For regular graphs, the new bound on the chromatic number is the same as the well-known Hoffman bound; however,
the two bounds are incomparable in general.

We prove that the eigenvalue bound for the max-$k$-cut is tight for several classes of graphs.
We investigate the presented bounds for specific classes of graphs, such as walk-regular graphs, strongly regular graphs,
and graphs from the Hamming association scheme.
\end{abstract}


\section{Introduction}

The max-$k$-cut problem is the problem of partitioning the vertex set of a graph into  $k$ subsets such that the total weight
of edges joining different sets is maximized.
The max-$k$-cut problem is also known as the minimum $k$-partition problem since minimizing the total weight of
the edges between vertices in the same part of the partition is equivalent to maximizing the $k$-cut.
The max-$k$-cut problem is NP-hard \cite{Arora}.
It has many applications such as VLSI design \cite{Barah88,Choetal98}, frequency planning \cite{Eisen02}, statistical physics \cite{Barah88},
digital-analogue convertors \cite{NesPolj86}, sports team scheduling \cite{Mitchell}, fault test generation \cite{Ibarra}, etc.

For the case that $k=2$, the max-$k$-cut problem is known as the max-cut problem.
The max-cut problem is one of the most studied combinatorial optimization problems and there is a large number of references related to the problem.
For studies on the cut polytope and its facets, see e.g., \cite{BarahMah,BoHamm91}.
In \cite{MoPo:09}, Mohar and Poljak  derived an eigenvalue bound for the max-cut problem.
A well-known semidefinite programming  (SDP) relaxation for the max-cut problem, in dual form,  was introduced by Delorme and Poljak \cite{DelPo93-2}.
The primal version of this basic SDP relaxation, one can find in e.g., \cite{PoljRendl952,GoemWill}.
In \cite{GoemWill}, Goemans and Williamson  showed  that  the basic SDP  relaxation for the max-cut has an error of no more than $13.82\%$.
Rendl, Rinaldi, and Wiegele \cite{ReRiWie10} incorporated this relaxation with additional inequalities within a branch and bound algorithm to solve --- to optimality --- max-cut instances of graphs with up to 100 vertices.

The max-$k$-cut problem was studied by Chopra and Rao \cite{ChopraRao93,ChopraRao95}, who derived several valid inequalities and facets for the $k$-partition polytope.
Further results on this polytope can be found in  \cite{DeGroLaus91,DeGroLa92,DezaLaur}.
To the best of our knowledge, there is no literature on eigenvalue bounds for the max-$k$-cut when $k\geq3$.
A SDP relaxation for the  max-$k$-cut problem was introduced by Frieze and Jerrum \cite{FriJerrum}.
In the same paper, the authors derived a polynomial-time approximation algorithm.
De Klerk, Pasechnik, and Warners \cite{KlerkPasWar} improved the approximation guarantees from  \cite{FriJerrum} for small values of $k$.
In \cite{Eisen02}, Eisenbl\"atter used a SDP relaxation for the minimum $k$-partition problem to show that the SDP relaxation of the $k$-partition polytope is strong.
Ghaddar, Anjos, and Liers \cite{GhAnjLie:08} developed a branch-and-cut algorithm that is based on a SDP relaxation for the minimum $k$-partition problem.
They were able to solve to optimality dense instances with up to 60 vertices and some special instances with up to 100 vertices, and for different values of $k$.
Anjos et al.~\cite{AjGhHuLiWi} improved the algorithm from  \cite{GhAnjLie:08} and developed a more efficient solver. \\

\noindent {\bf Main results and outline.}
In Section \ref{sect:SDP}  we derive a SDP relaxation for the max-$k$-cut problem from the SDP relaxation for the general graph partition problem from \cite{Sot:14}.
We show that the derived SDP relaxation is equivalent to the well-known SDP relaxation for the max-$k$-cut  problem  by Frieze and Jerrum \cite{FriJerrum}.
In Section \ref{sec:eigenbnd1} we derive a new eigenvalue bound for the  max-$k$-cut.
This eigenvalue bound is the first known closed form bound when $k\geq 3$ that is applicable to any graph,
and for $k=2$ it equals the eigenvalue bound for the max-cut by Mohar and Poljak  \cite{MoPo:09}.
The new eigenvalue bound is the optimal solution of a SDP relaxation that is dominated by the SDP relaxation from Section \ref{sect:SDP}.
In Section \ref{sect:eingBnd2} we improve the eigenvalue bound from Section \ref{sec:eigenbnd1}  by performing diagonal perturbations of the cost matrix.
This approach results in a new eigenvalue-based bound that can be obtained by solving a SDP problem.

Using the eigenvalue bound for the max-$k$-cut from Section \ref{sec:eigenbnd1},
we present in Section \ref{sec:chromatic} a bound on the chromatic number of a graph.
This is the first known closed form bound on the chromatic number of a graph that considers only eigenvalues of the Laplacian matrix of the graph.
For regular graphs, the new bound  equals the well-known  Hoffman bound \cite{Hoffman} on the chromatic number.
In general, however, there exist graphs for which our bound outperforms the Hoffman bound, and vice versa.

In Section \ref{sec:example} we demonstrate the quality of the presented bounds for several families of graphs, such as walk-regular graphs (see Section \ref{ExampleAssoc}),
strongly regular graphs (see Section \ref{Exam:SRG}),
and Hamming graphs (see Section \ref{ex:hamming}).
To our great surprise, we could not find in the literature any results on specific classes of graphs and their relation to the bounds for the max-$k$-cut with $k>2$.
We prove that the eigenvalue bound for the max-$k$-cut from Section \ref{sec:eigenbnd1} is tight for certain complete graphs, complete multipartite graphs
with $k$ color classes  of the same size, and certain graphs in the Hamming association scheme.
 The latter generalizes  a  result on the max-cut problem for binary `Hamming scheme graphs' by Alon and Sudakov \cite{AlonSudakov}.
 We also show that for certain $q$-ary Hamming scheme graphs and the max-$k$-cut  problem with $k\leq q$, the eigenvalue bound  from  Section \ref{sec:eigenbnd1} equals the SDP bound  from Section \ref{sect:SDP}.
 For walk-regular graphs, such as vertex-transitive graphs and graphs from association schemes,
 we prove that the two eigenvalue-based bounds for the  max-$k$-cut  from Section \ref{sec:eigenbnd1} and Section \ref{sect:eingBnd2} are equal.
 We also show that for walk-regular graphs  both eigenvalue bounds equal the optimal value of the SDP relaxation from Section \ref{sect:SDP} when $k=2$.
For strongly regular graphs, we also present a closed form expression for the optimal objective value of the SDP relaxation from Section \ref{sect:SDP}.
 This result is an extension of the results by De Klerk et al. \cite{deKlSotNaTr:10,dKPaDoSo:10} for the equipartition problem, and
 Van Dam and Sotirov \cite{DamSot13} for the general graph partition problem.
 We also show that for all strongly regular graphs except for the pentagon, the SDP bound from Section \ref{sect:SDP}  does not improve by adding triangle inequalities.

\vspace{1cm}

\section{SDP relaxations for the max-$k$-cut problem} \label{sect:SDP}

For a given undirected graph $G=(V,E)$, the max-$k$-cut problem asks for a partition of the vertex set $V$ into at most
$k$ subsets $V_i \subseteq V$ ($i=1,\ldots,k$) such that the weighted sum  of edges joining different sets is maximized.

Let $|V|=n$, let $A$ denote the adjacency matrix of $G$, and let $L=\Diag(Au_n)-A$ be its Laplacian matrix.
The well-known trace formulation of the max-$k$-cut problem is given by
\[
\begin{array}{rl}
\max & \frac{1}{2} \tr(X^{\mathrm T}LX)\\[1ex]
{\rm s.t.} & Xu_k=u_n \\[1ex]
&  x_{ij}\in \{0,1\}, ~~ i=1,\ldots, n,~~j=1,\ldots,k.
\end{array}
\]
Here $u_k$ and $u_n$ denote all-ones vectors of sizes $k$ and $n$, respectively.
Note that the columns of the matrix $X$ in this trace formulation are the incidence vectors of $V_i$ ($i=1,\ldots,k$).

In \cite{Sot:14} we derived a semidefinite programming relaxation for the
general graph partition problem. The general graph partition problem is the problem of partitioning vertices
of a graph into $k$ disjoint subsets of specified sizes such that the total weight of edges joining different sets  is optimized.
Thus,  a  mathematical model for the general graph partition problem reduces to a  mathematical model for the max-$k$-cut problem
after removing the constraints that impose restrictions on the subset sizes.

From the SDP relaxation for the general graph partition problem \cite{Sot:14}, we thus obtain the following SDP
relaxation for the max-$k$-cut problem (see also \cite{Rendl12}):
\begin{equation} \label{mainSDP}
\begin{array}{rl}
\max & \frac{1}{2} \tr(LY)\\[1ex]
{\rm s.t.} & \diag(Y) = u_n\\[1ex]
& kY-J_n\succeq 0, ~~Y\geq 0,
\end{array}
\end{equation}
where $J_n$ denotes the $n \times n$ all-ones matrix, and the `diag' operator maps a $n\times n$ matrix to the $n$-vector given by its diagonal.
Note that the SDP relaxation \eqref{mainSDP} contains a positive semidefiniteness constraint that is stronger than $Y\succeq 0$.
For the proof that $kY-J_n\succeq 0$ is a valid constraint, see e.g., \cite[Prop.~1]{Sot:14}.
We remind the reader that \eqref{mainSDP} can be obtained from the above trace formulation by introducing the variable $Y=XX^{\mathrm T}$, as usual.

For $k=2$, the nonnegativity constraints on the matrix variable $Y$ are implied by the positive semidefiniteness constraint and  $\diag(Y) = u_n$, see \cite{DamSot13}.
By using this fact we arrive at the following result.
\begin{prop}
For $k=2$, the SDP relaxation \eqref{mainSDP} is equivalent to the following well-known SDP relaxation for the max-cut problem by Delorme and Poljak \cite{DelPo93-2}:
\begin{equation}\label{FJ}
\begin{array}{rl}
\max & \frac{1}{4} \tr(LY)\\[1ex]
{\rm s.t.} & \diag(Y) = u_n\\[1ex]
& Y\succeq 0.
\end{array}
\end{equation}
\end{prop}
\noindent
{\em Proof}.
It is clear that there is a one-to-one correspondence between feasible $Y$ for \eqref{mainSDP} and feasible $Z$ for \eqref{FJ} by the relation $Z:=2Y-J_n$.
It is also easy to see that the two objectives coincide because $LJ_n=0$. \qed\\

Undoubtedly, the most cited SDP relaxation for the max-$k$-cut problem is the following relaxation by Frieze and Jerrum \cite{FriJerrum}:
\[
{\rm (FJ)}~~\begin{array}{rl}
\max & \frac{k-1}{2k} \tr(LY)\\[1ex]
{\rm s.t.} & \diag(Y) = u_n\\[1ex]
& Y\succeq 0, ~~Y \geq -\frac{1}{k-1}J_n.
\end{array}
\]
This SDP relaxation  is moreover exploited in \cite{FriJerrum} in a rounding heuristic to obtain a feasible solution of the max-$k$-cut problem.
It is not difficult to prove the following result.
\begin{prop}
The SDP relaxations \eqref{mainSDP}  and {\rm (FJ)} are equivalent.
\end{prop}

\noindent
{\em Proof}. The proof follows again by direct verification after applying appropriate variable transformations, see also \cite[Thm.~4]{Sot:14}. \qed\\

The SDP relaxation  \eqref{mainSDP} might be tightened by adding valid inequalities. For instance one can add triangle inequalities of type
\begin{equation}\label{triangle}
y_{ij}+y_{ik} \leq 1+y_{jk}, ~~\forall (i,j,k).
\end{equation}
For a given triple $(i,j,k)$ of (distinct) vertices, the constraint \eqref{triangle} ensures that
if $i$ and $j$ are in the same set of the partition and so are $i$ and $k$, then also $j$ and $k$ have to be in the same set.
The above mentioned inequalities are facet defining inequalities of the boolean quadric polytope, see e.g., \cite{Pad:89}.

One can also add  the independent set constraints
\begin{equation} \label{indepSet}
\sum\limits_{i<j, ~i,j\in Q} y_{ij}\geq 1, ~\mbox{for all}~ Q ~\mbox{with}~ |Q|=k+1.
\end{equation}
These constraints ensure that the graph with adjacency matrix $Y=(y_{ij})$ has no independent set ($Q$) of size $k+1$.
Anjos et al.~\cite{AjGhHuLiWi} use a bundle method to solve (FJ)  with  triangle and  independent set inequalities
in each node of  a branch-and-bound tree. Their approach resulted in a significantly faster  max-$k$-cut solver than the one presented in   \cite{GhAnjLie:08}.

Van Dam and Sotirov \cite{DamSot13} showed  how to aggregate triangle and independent set constraints  for graphs with symmetry in the context of a SDP relaxation for
the general graph partition problem.  In particular, their numerical results  show that for highly symmetric graphs one can compute a SDP bound that
includes all triangle and/or independent set inequalities  for graphs on about 100 vertices in a few seconds.
A similar approach can be applied here, see Section \ref{sec:example}.

\section{New bounds} \label{sect:eigen}

In this section, we derive two eigenvalue-based bounds for the max-$k$-cut and a new bound on the chromatic number of a graph.
In Section \ref{sec:eigenbnd1}, we derive an eigenvalue bound for the max-$k$-cut as the optimal solution of a SDP relaxation that is weaker than  \eqref{mainSDP}, in general.
Then, we  investigate possible improvements of the mentioned eigenvalue bound by perturbing the diagonal of the Laplacian matrix, see Section \ref{sect:eingBnd2}.
We show that such perturbations lead to a stronger bound.
In Section \ref{sec:chromatic} we derive a new lower bound on the chromatic number of a graph using the eigenvalue bound from Section \ref{sec:eigenbnd1}
for the max-$k$-cut.

\subsection{A new eigenvalue bound for the max-$k$-cut} \label{sec:eigenbnd1}

To derive a new  eigenvalue bound for the max-$k$-cut, we use a similar approach as the one in \cite{DamSot13}, where a matrix $*$-algebra called  the Laplacian algebra ${\mathcal L}$ is introduced.
A matrix $*$-algebra is a set of matrices that is closed under addition, scalar multiplication, matrix multiplication, and taking conjugate transposes.
It is well known that one can restrict optimization of a SDP problem to feasible points in a matrix $*$-algebra
that contains the data matrices of that problem as well as the identity matrix, see e.g.,  \cite{{Schrij79},{GoeRend:99},GatPa:04,deKlerk12}.

The  Laplacian algebra has a basis of  matrices that are obtained  from an orthonormal basis of eigenvectors corresponding to the eigenvalues of the Laplacian matrix $L$.
In particular, let $0=\lambda_0 \leq \lambda_1 < \ldots < \lambda_m=: \lambda_{\max}(L)$ be the distinct eigenvalues of $L$, where we allow that $\lambda_0 = \lambda_1$ in order to cater for the case that the graph is not connected.
Also, let $U_i$ be a matrix whose columns form an orthonormal basis of the eigenspace corresponding to the $i$th eigenvalue  $\lambda_i$ and
$F_i=U_iU_i^{\mathrm T}$ for $i=0,\ldots, m$, except in the case that the graph is not connected, in which case we split the eigenspace corresponding to eigenvalue $0$. In particular, we take $U_0=\frac 1{\sqrt{n}}u_n$ and $U_1$ a matrix whose columns form an orthonormal basis of the orthogonal complement of $\langle u_n \rangle$ within the eigenspace of eigenvalue $0$. Moreover, let $f_i =\rk F_i$ be the corresponding multiplicities, for $i=0,\ldots, m$.
The Laplacian algebra $\mathcal L$ is now defined as the span of $\{F_0, \dots, F_m\}$. This basis is called the basis of idempotents of $\mathcal L$, and  satisfies the following properties:
\begin{itemize}
\item $\sum\limits_{i=0}^m F_i=I$, $\sum\limits_{i=0}^m \lambda_iF_i~=L$
\item $F_iF_j = \delta_{ij} F_i$, ~$\forall i,j$
\item $F_i=F_i^*$, ~$\forall i$
\item $\tr(F_i)=f_i$, ~$\forall i$
\item $F_0= \frac{1}{n} J_n$.
\end{itemize}

Since the SDP relaxation \eqref{mainSDP} cannot be restricted to feasible points in ${\mathcal L}$
(for example, $\diag(Y)=u_n$ requires that the matrices in ${\mathcal L}$ have constant diagonal; we will indeed require this in Section \ref{ExampleAssoc}),
we relax several constraints in it.
In particular, we relax $\diag(Y)=u_n$ to $\tr(Y)=n$ and remove nonnegativity constraints, which leads to the relaxation:
\begin{equation} \label{eigSDP}
\begin{array}{rl}
\max & \frac{1}{2} \tr(LY)\\[1ex]
{\rm s.t.} & \tr(Y) = n\\[1ex]
& kY-J_n\succeq 0.
\end{array}
\end{equation}
Now, by restricting optimization of \eqref{eigSDP} to feasible points in  ${\mathcal L}$ we obtain the following result.
\begin{thm} \label{thmEig}
Let $G$ be a graph on $n$ vertices and $k$ an integer such that $2 \leq k\leq n$.
Then the optimal value of the  SDP relaxation  \eqref{eigSDP} is equal to
\begin{equation} \label{eigenvBnd}
\frac{n(k-1)}{2k} \lambda_{\max}(L),
\end{equation}
where $\lambda_{\max}(L)$ is the largest eigenvalue of  the Laplacian matrix $L$ of $G$.
\end{thm}

\noindent {\em Proof.} Because the Laplacian algebra ${\mathcal L}$ contains $L,I$, and $J_n$, there exists an optimal solution $Y$ to  \eqref{eigSDP} in ${\mathcal L}$ (see  e.g., \cite{{Schrij79},deKlerk12}), so
we may assume  $Y=\sum_{i=0}^m y_iF_i$ where $y_i \in \R$ ($i=0,\ldots, m$).
By exploiting the fact that $L=\sum_{i=0}^m \lambda_iF_i$ we have
\begin{equation} \label{obj}
\frac{1}{2}\tr(LY)= \frac{1}{2} \sum_{i=0}^m \lambda_i f_iy_i.
\end{equation}
Continuing in the same vein, from  $\tr(Y) = \sum_{i=0}^m y_i f_i = n$  and $kY-J_n\succeq 0$ we have
\begin{equation} \label{constr}
k-1 - \frac{k}{n}\sum_{i=1}^m f_iy_i\geq 0 \quad {\rm and} \quad y_i\geq 0 ~~(i=1,\ldots,m).
\end{equation}
In other words, our SDP relaxation \eqref{eigSDP} reduces to a linear programming (LP) problem with objective function  \eqref{obj} and constraints \eqref{constr}.
Now, it is not difficult to see that the optimal value of the resulting LP problem is $n(k-1)\lambda_{\max}(L)/2k$. \hfill\qed \\

We note that although our results also apply to graphs that are not connected, it is typically better to apply these results to the connected components of such graphs separately. This also applies to the eigenvalue bound \eqref{eigenvBnd} because the largest Laplacian eigenvalue of a graph is the maximum of the largest Laplacian eigenvalues of its connected components.

There are several interesting things related to the spectral  bound of Theorem \ref{thmEig}.
First, this is, to the best of our knowledge, the first closed form eigenvalue bound for the max-$k$-cut when $k>2$.
Second, for $k=2$ our eigenvalue bound coincides with the well-known eigenvalue bound for the max-cut by Mohar and Poljak \cite{MoPo:09}.
Our bound is particularly interesting for other small values of $k$.

Donath and Hoffman \cite{DonathHoffman} derived an eigenvalue-based bound for the general graph partition problem
that includes $k$  eigenvalues of a (perturbed) Laplacian matrix, where $k$ is the number of partition sets.
Their eigenvalue bound can also be reformulated  as an SDP problem, see  Alizadeh \cite{Alizadeh}.
From  the  Donath--Hoffman bound  one can derive an eigenvalue bound for the max-$k$-cut by bounding  all $k$ eigenvalues by the largest one.
This results in the bound $\frac{n}{2}\lambda_{\max}(L)$, which is weaker than the bound of Theorem \ref{thmEig}.
Further, it is known that for the equipartition problem  the eigenvalue bound from \cite{DonathHoffman}  is dominated by \eqref{mainSDP}, see \cite{KarRend:98,Sot:14}.
In Section \ref{sec:example} we provide several graphs for which our eigenvalue bound is tight.

\subsection{Strengthening the eigenvalue bound } \label{sect:eingBnd2}

In the following  we investigate possible improvements of the eigenvalue bound \eqref{eigenvBnd}.
In particular, we show that {certain perturbations of the diagonal of the Laplacian matrix
 lead to a bound that is in general stronger than \eqref{eigenvBnd} but weaker than \eqref{mainSDP}.
 Diagonal perturbations of eigenvalue optimization problems and their relations to semidefinite programming  were  studied in  \cite{Alizadeh}.
Perturbations of the cost matrix by a diagonal matrix with zero trace do not
change the optimal value of the max-$k$-cut problem, but have an effect on the maximal eigenvalue of the Laplacian matrix.
This was pointed out already by  Delorme  and Poljak \cite{DelPo93} for  the max-cut problem.
Therefore, we investigate the following optimization problem:
\begin{equation} \label{imprEig}
\min_{d^{\mathrm T}u_n=0}\frac{n(k-1)}{2k} \lambda_{\max}(L + \Diag(d)),
\end{equation}
where `Diag' is the adjoint operator of `diag'.
The vector $d$ is known as the correcting vector, see e.g., \cite{DelPo93}.
The eigenvalue optimization problem \eqref{imprEig} can be formulated as  the following semidefinite program:
\begin{equation*} \label{DualimprEigSDP}
\begin{array}{rl}
\min & \frac{n(k-1)}{2k}\mu \\[1ex]
{\rm s.t.} & \mu I_n - (L + \Diag(d))=Z\\[1ex]
& Z \succeq 0, ~~   d^{\mathrm T}u_n=0,
\end{array}
\end{equation*}
whose dual problem is
\begin{equation} \label{imprEigSDP}
\begin{array}{rl}
\max & \frac{1}{2} \tr(LY) \\
{\rm s.t.}  & \diag(Y)=\frac{k-1}{k} u_n \\
& Y\succeq 0.
\end{array}
\end{equation}
Note that  the optimal value of the eigenvalue problem \eqref{imprEig} equals  the optimal value of the SDP relaxation \eqref{imprEigSDP}.
In what follows, we  show that  the SDP relaxation \eqref{mainSDP} dominates \eqref{imprEigSDP}.
By the variable transformation $Z=Y+ \frac{1}{k}J_n$ together with $LJ_n=0$, we obtain from \eqref{imprEigSDP}  the following equivalent  SDP relaxation:
\[
\begin{array}{rl}
\max & \frac{1}{2} \tr(LZ)\\[1ex]
{\rm s.t.} & \diag(Z) = u_n\\[1ex]
& kZ-J_n\succeq 0.
\end{array}
\]
It is now clear that the SDP relaxation  \eqref{mainSDP} dominates \eqref{imprEigSDP} because of the constraints $Z \geq 0$ (and confirm also that \eqref{imprEigSDP} dominates \eqref{eigSDP}).

\begin{prop} \label{1to11}
 The SDP relaxation \eqref{mainSDP} is equivalent to the   SDP relaxation  \eqref{imprEigSDP} with  additional inequalities $Y \geq -\frac{1}{k} J_n$.
\end{prop}

We summarize the relations between the presented SDP relaxations of the max-$k$-cut problem and corresponding bounds as follows:
\[
\begin{array}{ccccccc}
\mbox{SDP }  \eqref{eigSDP} & \geq & \mbox{SDP } \eqref{imprEigSDP}   & \geq &  \mbox{SDP } \eqref{mainSDP}  &  \geq  & \mbox{max-$k$-cut} \\
\| && \| &&  &&\\
\mbox{$\lambda_{\max}$-bound } \eqref{eigenvBnd} && \mbox{Perturbed `bound' }  \eqref{imprEig}  &&
\end{array}
\]
where $A\geq B$ means that $B$ dominates $A$.
For $k=2$, the nonnegativity constraints in  \eqref{mainSDP} are redundant (see \cite{DamSot13}),
and  therefore the SDP relaxations \eqref{imprEigSDP} and \eqref{mainSDP} are equivalent.
In Section \ref{ExampleAssoc}, we prove that for walk-regular graphs, the SDP bounds \eqref{eigSDP} and  \eqref{imprEigSDP} are equal.
In Section \ref{ex:hamming}, we prove that for certain graphs in the $q$-ary Hamming association scheme with $k \leq q$, the SDP bounds \eqref{eigSDP} and \eqref{mainSDP} are equal,
while on top, the SDP bound \eqref{eigSDP} is tight for $k=q$.

\subsection{A new bound on the chromatic number} \label{sec:chromatic}

We will now derive a bound on the chromatic number of a graph  from the eigenvalue bound \eqref{eigenvBnd} for the max-$k$-cut problem.
We show that for regular graphs our bound equals  the well-known Hoffman bound \cite{Hoffman}.

A coloring of a graph is an assignment of colors to the vertices of the graph such that no two adjacent vertices have  the same color.
The smallest number of colors needed to color a graph $G$ is called its chromatic number $\chi(G)$.
Since a coloring with $k$ colors  is the same as a partition of the vertex set into $k$ independent sets, we have the following observation.
For a given graph $G=(V,E)$ and integer $k$,
\begin{center}
if ~ max-$k$-cut $< |E|$, ~then~ $\chi(G) \geq k+1$.
\end{center}
By using this in combination with the eigenvalue bound \eqref{eigenvBnd}, we obtain the following bound on the chromatic number of $G$:

\begin{thm} \label{thm:Chromatic} Let $G=(V,E)$ be a graph with Laplacian matrix $L$. Then
\begin{equation} \label{Bnd:color}
\chi(G) \geq  1 + \frac{2|E|}{n \lambda_{\max}(L) - 2|E| }.
\end{equation}
\end{thm}

\noindent {\em Proof.} Let $k=\lceil \frac{2|E|}{n \lambda_{\max}(L) - 2|E| }\rceil$, then $k<1+\frac{2|E|}{n \lambda_{\max}(L) - 2|E| }$. From this, it follows that the upper bound \eqref{eigenvBnd} for the max-$k$-cut is less than $|E|$, and hence \eqref{Bnd:color} follows.
\qed \\

For a $\kappa$-regular graph with $n$ vertices and adjacency matrix $A$ with smallest eigenvalue $\theta_{\min}$,
we have that $\kappa = 2|E|/n$ and $\lambda_{\max}(L)=\kappa - \theta_{\min}$ because $L=\kappa I-A$.
Thus, for a $\kappa$-regular graph, \eqref{Bnd:color} can be rewritten as
\[
\chi(G) \geq 1 - \frac{\kappa}{\theta_{\min} },
\]
which is exactly the well-known Hoffman bound \cite{Hoffman} on the chromatic number (when restricted to regular graphs).
The new bound can be strictly better than the Hoffman bound.
For instance for the complete graph on $100$ vertices minus an edge, our bound is 99 (which is also the chromatic number) while the Hoffman bound is $51$.

Recently, several generalizations  of the Hoffman bound were presented that include several or all eigenvalues of the adjacency matrix
and/or Laplacian matrix, see \cite{ElpWo:14}. Our new bound differs from these since it includes only one eigenvalue of the Laplacian matrix.

\section{Specific classes of graphs and the max-$k$-cut} \label{sec:example}

In this section we investigate the introduced eigenvalue and SDP bounds for the max-$k$-cut for specific classes of graphs.
In particular, we show  that the eigenvalue bound \eqref{eigenvBnd}
is tight for some complete graphs, complete multipartite graphs with $k$ color classes of the same size, and certain graphs in the $q$-ary Hamming association scheme with $k=q$.
We show that for the Coxeter graph  the SDP relaxation \eqref{mainSDP} with  $k=2$  and additional triangle and independent set inequalities is tight.
In Section \ref{ExampleAssoc}, we prove that for walk-regular graphs  the eigenvalue bound  \eqref{eigenvBnd}
is equal to the eigenvalue bound \eqref{imprEig}.
In Section \ref{Exam:SRG}, we derive a closed form bound for the max-$k$-cut of strongly regular graphs from the SDP relaxation \eqref{mainSDP},
and give conditions under which it is equal to the eigenvalue bound \eqref{eigenvBnd}.
In Section \ref{ex:hamming}, we consider graphs from the Hamming association scheme. We now start with an easy example.

\begin{example} {\em Complete multipartite graph with $k$ color classes of the same size}.\\
The largest Laplacian eigenvalue  of $K_{k\times m}$ is $k m$, and hence the eigenvalue bound \eqref{eigenvBnd} equals
\begin{equation*} \label{Kkm}
\frac{k(k-1) m^2}{2},
\end{equation*}
which is clearly the max-$k$-cut in $K_{k\times m}$ (and the number of edges).
\end{example}

\begin{example} {\em Complete graphs.} \\
The largest Laplacian eigenvalue of  the complete graph $K_n$  is $n$.
Thus, the eigenvalue upper bound for the max-$k$-cut problem  \eqref{eigenvBnd} is  $\frac{n^2(k-1)}{2k}$.
To derive the max-$k$-cut in $K_n$, we need to solve the following optimization problem:
\[
\begin{array}{rl}
\max & \sum\limits_{i=1}^k \sum\limits_{j=i+1}^k m_i m_j \\
{\rm s.t.} & \sum\limits_{i=1}^k m_i = n\\
& m_i \in \mathbb{Z}, ~m_i \geq 0, ~ i=1,\ldots,k.
\end{array}
\]
Under the given constraints, it is not hard to show that
$$\sum\limits_{i=1}^k \sum\limits_{j=i+1}^k m_i m_j=\frac{n^2(k-1)}{2k}-\frac 12 \sum_{i=1}^k(m_i-\frac nk)^2,$$
and hence that if $n$ is a multiple of $k$, then the max-$k$-cut equals the eigenvalue bound.

If $n$ is not a multiple of $k$, then the max-$k$-cut is attained by taking $m_i=\lceil \frac{n}{k}  \rceil$ ($n-k\lfloor n/k \rfloor$ times) and  $m_i=\lfloor \frac{n}{k} \rfloor$ ($k-n+k\lfloor n/k \rfloor$ times). We can in fact show (but omit the elaborate but straightforward details) that if $n \equiv e \mod k$ with $0 \leq e \leq k-1$, then the max-$k$-cut equals
$$\frac{n^2(k-1)}{2k}-\frac{e(k-e)}{2k}.$$
Thus, as long as $\frac{e(k-e)}{2k}<1$, the rounded eigenvalue bound equals the max-$k$-cut. This holds for example for $e=1$ or $2$ and all $k$, but also for $k \leq 7$ and all $n$. On the other hand, it shows that the gap between the eigenvalue bound and the max-$k$-cut can be arbitrary large (take $e=\lfloor k/2 \rfloor$ and let $k \rightarrow \infty$). The `smallest' case where the (rounded) eigenvalue bound is not the max-$k$-cut is the case that $n=12$ and $k=8$.
\end{example}

The following example shows that after adding the triangle inequalities  \eqref{triangle} and independent set constraints \eqref{indepSet} to \eqref{mainSDP}, the resulting bound can be tight.

\begin{example}{\em The Coxeter graph.} \\
The Coxeter graph is a regular graph with $28$ vertices and $42$ edges.
Delorme and Poljak \cite{DelPo93} proved that the optimal value of the max-cut for the Coxeter graph is $36$,
while the eigenvalue bound \eqref{eigenvBnd} is the same as the SDP bound \eqref{mainSDP}, which equals $37.89$.
For $k=2$, we computed \eqref{mainSDP} with additional triangle  inequalities  \eqref{triangle},  and also
with additional  triangle  and  independent set inequalities  \eqref{indepSet}. The obtained  bounds are  $36.75$ and $36$, respectively,
and hence they are tight (the first one after rounding, of course).

Note that we computed the SDP relaxations with additional inequalities by aggregating the constraints as described in \cite{DamSot13}.
The relaxation with all triangle and independent set inequalities was solved in only $0.14$ s.
\end{example}

\subsection{Walk-regular graphs and association schemes} \label{ExampleAssoc}

We show here that for so-called walk-regular graphs,
the eigenvalue bound \eqref{eigenvBnd} is equal to the improved eigenvalue bound \eqref{imprEig}.

A graph with adjacency matrix $A$ is called walk-regular if $A^{\ell}$ has constant diagonal for every nonnegative integer $\ell$.
The matrix $A^{\ell}$ contains the numbers of walks of length $\ell$ between vertices, so the definition is equivalent to requiring that
the number of walks of length $\ell$ from a vertex to itself is the same for every vertex. The class of walk-regular graphs contains
all vertex-transitive graphs, distance-regular graphs, and graphs whose Laplacian matrix is contained in the Bose-Mesner algebra of an association scheme,
among others (such as the graph from \cite[Fig.~2]{GM}). For more on walk-regular graphs, we refer to \cite{GM} or \cite{Gbook}.

Because the diagonal of $A^2$ contains the vertex degrees of the graph, it follows that a walk-regular graph is regular, say of degree $\kappa$.
Consider now the Laplacian algebra ${\mathcal L}$ of Section \ref{sect:eigen}.
Because a walk-regular graph is regular, its Laplacian matrix equals $L=\kappa I-A$.
Now it follows that all matrices in the Laplacian algebra ${\mathcal L}$ of a walk-regular graph have constant diagonal, in particular the
idempotent matrices $F_0,F_1,\dots, F_m$. Below we will exploit the fact that $F_m$ has constant diagonal.

Recall that $m$ is the index of the largest eigenvalue of the Laplacian matrix $L$, so that
\begin{equation} \label{LE}
L F_m = \lambda_{\max}(L)\cdot F_m.
\end{equation}
We now define
\[
\bar{Y} :=\frac{n(k-1)}{k \tr(F_m)} F_m,
\]
then $\bar{Y}$ has constant diagonal, with
\[
\diag(\bar{Y})=\frac{k-1}{k}u_n,
\]
and hence $\bar{Y}$ is feasible for \eqref{imprEigSDP}, with objective value
\[
\frac{1}{2} \tr(L\bar{Y})=\frac{n(k-1)}{2k \tr(F_m)} \tr(L F_m) = \frac{n(k-1)}{2k} \lambda_{\max}(L),
\]
where we exploited the fact that $L = \sum_{i=0}^m \lambda_i F_i$ and \eqref{LE}. Thus, the maximal value of the SDP bound \eqref{imprEigSDP} is attained in $\bar{Y}$, and \eqref{imprEigSDP} (and hence \eqref{imprEig}) equals the eigenvalue bound \eqref{eigenvBnd}. We remark that this implies also that the optimum correcting vector $d$ in \eqref{imprEig} is equal to the zero vector. Poljak and Rendl \cite{PoljRendl951} proved this result for the max-cut problem in vertex-transitive graphs.

Note that the SDP relaxation \eqref{mainSDP} is in general stronger than \eqref{imprEigSDP} because it contains additional inequality constraints.
Since for $k=2$ the nonnegativity constraints in \eqref{mainSDP} are redundant,
the eigenvalue bound \eqref{eigenvBnd} for the max-cut  is equal to the optimal solution of  \eqref{mainSDP} for all walk-regular graphs.
Goemans and Rendl \cite{GoeRend:99} proved this for the max-cut problem in graphs in an association scheme, and therefore our result generalizes theirs.
We summarize our results in the following theorem.
\begin{thm}
Let $G$ be a walk-regular graph on $n$ vertices and let $k$ be an integer such that $2\leq k \leq n$.
Then the eigenvalue bound \eqref{eigenvBnd} equals the improved eigenvalue bound \eqref{imprEig}.
Moreover, for $k=2$ the eigenvalue bound \eqref{eigenvBnd} equals the optimal value of the SDP relaxation \eqref{mainSDP}.
\end{thm}

What we required in the above analysis is that the idempotent $F_m$ has constant diagonal.
We remark that this can also be the case for graphs that are not walk-regular. In fact, this holds for every regular bipartite graph.
However, for bipartite graphs, the max-$k$-cut problem is trivial.

\subsection{Strongly regular graphs} \label{Exam:SRG}
Strongly regular graphs  belong to the class of walk-regular graphs analyzed in Section \ref{ExampleAssoc}.
Here we derive a closed form expression for the optimal objective value of the SDP relaxation \eqref{mainSDP}.
We also derive a condition under which the eigenvalue bound \eqref{eigenvBnd} equals the  optimal objective value of the SDP relaxation \eqref{mainSDP}.
A similar approach was used in \cite{deKlSotNaTr:10,dKPaDoSo:10} to derive an eigenvalue bound for the equipartition and
in \cite{DamSot13} to derive an eigenvalue bound for the general graph partition problem.

A  $\kappa$-regular graph $G=(V,E)$
on $n$ vertices is called strongly regular with parameters $(n, \kappa, \lambda, \mu)$ whenever it is
not complete or edgeless and every two distinct vertices have $\lambda$ or $\mu$ common neighbors, depending
on whether the two vertices are adjacent or not, respectively.
Let $A$ be the adjacency matrix of a strongly regular graph $G$ with parameters $(n, \kappa, \lambda, \mu)$.
Since $G$ is regular with valency $\kappa$, we have that $\kappa$ is an eigenvalue of $A$ with eigenvector $u_n$.
The matrix $A$ has exactly two distinct eigenvalues associated with eigenvectors orthogonal to $u_n$.
These two eigenvalues are known as the restricted eigenvalues and are usually denoted by $r$ and $s$, where $r \geq 0$ and $s<0$.
Hence, the largest Laplacian eigenvalue $\lambda_{\max}(L)$ of $G$ equals $\kappa-s$.

Now, by applying the general theory of symmetry reduction to the SDP problem \eqref{mainSDP}  as described in, e.g.,
\cite{DamSot13,deKlSotNaTr:10,dKPaDoSo:10} we obtain the following result.
\begin{thm} \label{strgClosedForm}
Let  $G=(V,E)$ be a  strongly regular graph with parameters $(n,\kappa,\lambda,\mu)$ and restricted eigenvalues $r \geq
0$ and $s<0$. Let $k$  be an integer such that $2 \leq k<n$. Then the SDP bound \eqref{mainSDP} for the max-$k$-cut of $G$ is given by
\[
\min \left \{ \frac{n(k-1)}{2k}(\kappa - s), ~\frac{1}{2}\kappa n \right \}.
\]
\end{thm}
Note that the objective value of \eqref{mainSDP} is equal to  the eigenvalue bound \eqref{eigenvBnd} unless  $\frac{k-1}{k} > \frac{\kappa}{\kappa -s}$.
We remark that the term $\frac{n(k-1)}{2k}(\kappa - s)$ is indeed the eigenvalue bound \eqref{eigenvBnd}, which by the results in the previous section
equals the SDP bound \eqref{imprEigSDP}. The term $\frac{1}{2}\kappa n$ in the above bound comes in because of the nonnegativity constraints in \eqref{mainSDP}.
Note also that this term equals the number of edges $|E|$, and hence is an obvious upper bound for the max-$k$-cut.

In \cite{DamSot13,GoeRend:99}, it is explained how one can aggregate the triangle inequalities \eqref{triangle} for graphs with symmetry.
Further, in \cite{DamSot13} it is shown that  in the case of the graph partition problem for strongly regular graphs, and
after simplifying and removing equivalent inequalities from the set of $3\binom{n}{3}$ triangle inequalities,  there remain
at most four inequalities. Clearly, the same inequalities remain for the max-$k$-cut problem.
Therefore, following  similar arguments as those in \cite{DamSot13}, we conclude that
adding triangle inequalities to \eqref{mainSDP} does {\em not} improve the bound for strongly regular graphs, except possibly for the pentagon.
Indeed, if we consider the pentagon then the SDP bound \eqref{mainSDP} with $k=2$ equals $4.52$, and after adding all triangle inequalities to \eqref{mainSDP} the bound becomes $4.16$.

We also tested the quality of the bound obtained after adding independent set inequalities, for the case $k=2$.
Note that it is not clear how to aggregate those constraints for $k>2$, see \cite{DamSot13}.
Our numerical results show that after adding constraints \eqref{indepSet} to \eqref{mainSDP}, the bound improves for some instances.
For example, for the Kneser graph $K(6,2)$ the eigenvalue bound \eqref{eigenvBnd} equals the SDP bound  \eqref{mainSDP}, which is equal to $33.75$,
while after adding aggregated independent set constraints to \eqref{mainSDP} the bound is $30$.
For the pentagon, adding all triangle and independent set constraints to the SDP relaxation \eqref{mainSDP} gives a bound of $4$.

\subsection{The Hamming graphs} \label{ex:hamming}

We now consider (the walk-regular) graphs from the Hamming association scheme. With vertices represented by $d$-tuples of letters from an alphabet of size $q$,
the adjacency matrices $H(d,q,j)$ ($j=0,\ldots, d$) of the Hamming association scheme are defined by the number of positions in which two $d$-tuples differ.
In particular, let $V$  be the set of all vectors of length $d$ over an alphabet  of size $q$, so $n=q^d$.
Then, $H(d,q,j)_{x,y}=1$ if $x$ and $y$ differ in $j$ positions, for $x,y\in V$. The graph $H(d,q,1)$ is the well-known Hamming graph, which can also
be obtained as the Cartesian (or direct) product of $d$ copies of the complete graph $K_q$.
The eigenvalues of the graphs in the Hamming scheme can be obtained using Kravchuk polynomials, see \cite[Thm.~30.1]{LintWilson}, i.e., $H(d,q,j)$ has eigenvalues
$$K_j(i)=\sum_{h=0}^j(-q)^h(q-1)^{j-h}{d-h \choose j-h}{i \choose h},$$
with multiplicities $m_i={d \choose i}(q-1)^i$, for $i=0,\dots,d$.
Note that $H(d,q,j)$ is a regular graph with degree $K_j(0)={d \choose j}(q-1)^j$, and hence its Laplacian eigenvalues are $K_j(0)-K_j(i)$, for $i=0,\dots,d$.

For the particular case that $j=d$, we have that $H(d,q,d)=(J_q-I_q)^{\otimes d}$, the $d$-th power of $J_q-I_q$ using the Kronecker product, which
confirms that the eigenvalues of $H(d,q,d)$ are $K_d(i)=(-1)^i(q-1)^{d-i}$, for $i=0,\dots,d$. Moreover, this particular form allows us to write the idempotent matrix of Section \ref{sec:eigenbnd1} that corresponds to the largest Laplacian eigenvalue $\lambda_{\max}(L)$ in a useful form, as we shall see later. Because
the graph is regular of degree $K_d(0)=(q-1)^d$ with smallest eigenvalue $\theta_{\min}=K_d(1)=-(q-1)^{d-1}$, it follows that
\begin{equation*} \label{lamMaxH}
\lambda_{\max}(L)=(q-1)^d + (q-1)^{d-1}=q(q-1)^{d-1}.
\end{equation*}

Alon and Sudakov \cite{AlonSudakov} showed that the optimal value for the max-cut problem for $H(d,2,j)$
coincides with the optimal solution of \eqref{mainSDP} (and implicitly, also with the eigenvalue bound \eqref{eigenvBnd}), when $j$ is even, $j>d/2$, and $d$ is large enough with $j/d$ fixed.
 We will extend this result to the case $q>2$. Note that for $q=2$, the graphs $H(d,2,j)$ are bipartite for $j$ odd, and so for these graphs, the max-cut problem is not relevant. For $q>2$, all graphs $H(d,q,j)$ are however non-bipartite and also connected.

Consider again $H(d,q,d)$. The idempotent matrices of Section \ref{sec:eigenbnd1} for $H(1,q,1)=J_q-I_q$ are $\frac1qJ_q$ (with eigenvalue $q-1$) and $I_q-\frac1qJ_q$ (with eigenvalue $-1$). Since $H(d,q,d)$ is obtained by taking Kronecker products of $H(1,q,1)$, its idempotents can be obtained using Kronecker products of these two `basic' idempotents. It thus follows that the idempotent for $\lambda_{\max}(L)$ of $H(d,q,d)$ equals
\begin{equation*} \label{FnmH}
F = \frac{1}{q^d} \sum_{i=0}^{d-1} J_q^{\otimes i}\otimes(q I_q -J_q)\otimes J_q^{\otimes d-1-i}=\frac{1}{q^d}(q \sum_{i=0}^{d-1} J_{q^i}\otimes I_q\otimes J_{q^{d-1-i}}-dJ_{q^d}),
\end{equation*}
with $\tr (F)=d(q-1)$.
This is in fact not just an idempotent for $H(d,q,d)$, but it is also an idempotent for the Hamming association scheme, and for all graphs $H(d,q,j)$, $j=1,\dots,d$.
The corresponding Laplacian eigenvalue $\lambda$ of $H(d,q,j)$ equals
\begin{equation} \label{lambda}
\lambda=K_j(0)-K_j(1)=q(q-1)^{j-1} {d-1 \choose j-1}.
\end{equation}
Note that in order to avoid confusion, we didn't use an index for $F$ and $\lambda$, because the ordering of the eigenvalues here is not increasing, as in Section \ref{sec:eigenbnd1}. We now make the following conjecture on $\lambda$.
\begin{conj}\label{conj:lambda} Let $q \geq 2$ and $j\geq d-\frac{d-1}q$, with $j$ even if $q=2$. Then $\lambda$ is the largest Laplacian eigenvalue of $H(d,q,j)$.
\end{conj}
Alon and Sudakov \cite[Prop.~3.2]{AlonSudakov} showed that $\lambda$ is indeed the largest Laplacian eigenvalue of $H(d,2,j)$ for $j$ even, $j>d/2$, when $d$ is large enough with $j/d$ fixed.
We confirmed the conjecture numerically (comparing similarly as Alon and Sudakov the values of the Kravchuk polynomials, i.e., that $K_j(i) \geq K_j(1)$ for all $i=0,\dots,d$) for all pairs $(d,q)$ with $d\leq 30$ and $q\leq15$. We have no hard proof however, except for the case $j=d$, of course.

Under the condition that $\lambda$ is indeed the largest Laplacian eigenvalue of $H(d,q,j)$, this would imply that the eigenvalue bound \eqref{eigenvBnd} for the max-$k$-cut equals
$$\frac{q^d(k-1)}{2k}\lambda=\frac{k-1}{2k}q^{d+1}(q-1)^{j-1} {d-1 \choose j-1}.$$
In Section \ref{ExampleAssoc}, we showed that the matrix
\[
\bar{Y}=  q^d \frac{k-1}{kd(q-1)} F,
\]
is feasible for \eqref{imprEigSDP}. Moreover, if $k\leq q$, then
\[
\bar{Y} \geq -\frac{k-1}{k(q-1)} J_{q^d}  \geq -\frac{1}{k}J_{q^d}.
\]
From Proposition \ref{1to11}, we thus obtain the following result.

\begin{thm}\label{thm:hammingsamebounds}
Let $k\leq q$, $j\geq d-\frac{d-1}q$, with $j$ even if $q=2$, and consider the graph $H(d,q,j)$. If $\lambda$ in \eqref{lambda} is the maximal Laplacian eigenvalue, then for the max-$k$-cut problem, the eigenvalue bound \eqref{eigenvBnd} and the SDP bound \eqref{mainSDP} are equal.
\end{thm}

Next, we will construct optimal $q$-cuts, showing that the eigenvalue bound \eqref{eigenvBnd} is tight for $k=q$. Because the adjacency matrix of  $H(d,q,d)$ is of the form
$(J_q - I_q) \otimes (J_q-I_q)^{\otimes d-1} $,
it follows that  the max-$q$-cut of this graph is equal to the number of edges
$\frac{1}{2} (q-1)^d q^d$,
which indeed proves that the eigenvalue bound is tight.
It also shows that the chromatic number is at most $q$. Together with the eigenvalue (lower) bound for the chromatic number,
this shows that the chromatic number of $H(d,q,d)$ equals $q$.

For the case $k=q=2$, Alon and Sudakov construct a max-cut using a specific eigenvector corresponding to the largest Laplacian eigenvalue (or to the smallest eigenvalue of the adjacency matrix, to speak in their terms). For $q=2$, the specific eigenvectors have entries only $1$ and $-1$, and these values induce the cut.
Looking at these eigenvectors more carefully, we can describe such a cut, and more generally a $q$-cut for $H(d,q,j)$ in a combinatorial way as follows. Recall that we can write the vertex set as $V=\{1,\dots,q\}^n$. Fix a coordinate, say the first, and partition $V$ in $q$ parts according to the value of the vertex in this coordinate, i.e., let $V_i=\{x \in V: x_1=i\}$, for $i=1,\dots,q$. It is quite straightforward to count the number of edges in this $q$-cut for $H(d,q,j)$ as follows. Let $i=1,\dots,q$, and consider a vertex $x$ in $V_i$. Let $h \neq i$. Then every vertex in $V_h$ that is adjacent to $x$ differs from $x$ in the first coordinate, so it should differ in $j-1$ of the remaining $d-1$ coordinates. This implies that $x$ has ${d-1 \choose j-1}(q-1)^{j-1}$ neighbors in $V_h$, and hence there are $q^{d-1}{d-1 \choose j-1}(q-1)^{j-1}$ edges between $V_i$ and $V_h$. The total number of edges in the $q$-cut is therefore
$$\frac12q(q-1)q^{d-1}{d-1 \choose j-1}(q-1)^{j-1}=\frac{n(q-1)}{2q}\lambda,$$
and hence the eigenvalue bound \eqref{eigenvBnd} is tight, again provided that $\lambda$ is indeed the largest Laplacian eigenvalue.

\begin{thm}\label{thm:hammingtight}
Let $j\geq d-\frac{d-1}q$, with $j$ even if $q=2$, and consider the graph $H(d,q,j)$. If $\lambda$ in \eqref{lambda} is the maximal Laplacian eigenvalue, then the optimal value of the max-$q$-cut equals the eigenvalue bound \eqref{eigenvBnd}.
\end{thm}

\section{Conclusion}

In this paper, we  presented several new bounds for the max-$k$-cut problem and the chromatic number of a graph,
 and analyzed their quality on specific classes of graphs.

In particular, in Theorem \ref{thmEig} we derived an eigenvalue bound for the max-$k$-cut problem for $k\geq 2$,
extending a  well-known eigenvalue bound for the max-cut problem by Mohar and Poljak \cite{MoPo:09}.
We exploited this new eigenvalue bound \eqref{eigenvBnd}
to derive a new lower bound  \eqref{Bnd:color} on the chromatic number of a graph.
This new bound, when applied to regular graphs, reduces to the well-known Hoffman bound \cite{Hoffman} on the chromatic number.
In general however, our bound is not dominated by the Hoffman bound, nor the other way around.
In Section \ref{sect:eingBnd2} we showed how to improve the eigenvalue bound  \eqref{eigenvBnd}. The resulting bound does not have a closed form
expression, but equals  the optimal solution of the SDP relaxation \eqref{imprEigSDP}.

For walk-regular graphs, we showed that the eigenvalue bound  \eqref{eigenvBnd} is equal to the improved eigenvalue bound \eqref{imprEig}.
For strongly regular graphs, we derived a closed form  bound that dominates the previously mentioned eigenvalue bounds.
Finally, under some assumption on the behavior of the Kravchuk polynomials (which we tested extensively, see Conjecture \ref{conj:lambda}), we showed that for $j\geq d-\frac{d-1}q$,  the optimal value of the max-$q$-cut equals the eigenvalue bound \eqref{eigenvBnd} for the graphs $H(d,q,j)$ in the Hamming association scheme.
This extends a similar result for the binary Hamming graphs and the max-cut problem by Alon and Sudakov \cite{AlonSudakov}.
Under the same assumption, we also proved that for the max-$k$-cut problem for $H(d,q,j)$, with $j\geq d-\frac{d-1}q$ and $k\leq q$, the eigenvalue bound \eqref{eigenvBnd} and the SDP bound \eqref{mainSDP} are equal.

\end{document}